\documentclass[a4paper,11pt]{article}
\usepackage{xy}
\usepackage{graphicx,amsmath,amssymb,amsthm,color}

\newtheorem{thm}{Theorem}[section]
\newtheorem{lemma}[thm]{Lemma}

\begin{document}

\begin{center}

{\huge \bf Homotopy types of truncated projective resolutions}

{W.H.Mannan}

\end{center}

\noindent {\bf Published:} 

\noindent Homology, Homotopy and Applications Vol. 9 (2007), No. 2, pp. 445-449

\bigskip
{MSC: 16E05.} \hfill Keywords: projective resolution, homotopy type

\bigskip
\begin{tiny}
{\bf Abstract} We work over an arbitrary ring $R$.  Given two truncated
projective resolutions of equal length for the same module we
consider their underlying chain complexes.  We show they may be
stabilized by projective modules to obtain a pair of complexes of
the same homotopy type.
\end{tiny}

\section{Introduction}

Truncated projective resolutions are of interest in both algebraic
geometry and algebraic topology.  If the modules in a resolution
of length $n$ are assumed to be free, then the $n^{\rm th}$
homology group is the $n^{\rm th}$ syzygy of the module being
resolved. The minimal possible dimensions of the modules in such
resolutions were of interest to mathematicians such as Hilbert and
Milnor (see \cite{Miln}).

In algebraic topology, truncated projective resolutions arise as the
algebraic complexes associated to $(n-1)$-connected universal covers
of CW-complexes of dimension $n$.  Of particular interest is the
case $n=2$, as classification of the homotopy types of these
truncated resolutions is closely related to Wall's D2 problem (see
the introduction to \cite{John}).

Given two truncated projective resolutions of the same module (of
equal  length), their final modules may be stabilized to produce
homotopy equivalent algebraic complexes.  This is a generalization
of Schanuel's lemma which merely equates the final homology
groups. The work of mathematicians such as Milnor, Whitehead and
Wall suggest they were familiar with this basic homological
result. Indeed, Wall's obstruction is suggestive of the modules
required to stabilize the complexes (see \cite{Wall}, \S3).

Given one truncated projective resolution of a module this result
provides a handle on all other possible truncated projective
resolutions of the same length.  Our purpose in this paper is to
provide a simple proof of the result by explicitly constructing
the desired homotopy equivalence between the two stabilized
algebraic complexes.

\bigskip
Formally, let $R$ be a ring with identity and let $M$ be a module
over $R$. We assume a right action on all modules.  Suppose we
have exact sequences:
\begin{center}$
P_n \stackrel{\partial_n}{\longrightarrow} P_{n-1}
\stackrel{\partial_{n-1}}{\longrightarrow}\cdots
\stackrel{\partial_2}{\longrightarrow} P_1 \stackrel{\partial_1}
{\longrightarrow} P_0 \stackrel{\epsilon}{\dashrightarrow} M
\dashrightarrow 0 $\end{center} and
\begin{center}$
Q_n \stackrel{\partial_n'}{\longrightarrow} Q_{n-1}
\stackrel{\partial_{n-1}'}{\longrightarrow} \cdots
\stackrel{\partial_2'}{\longrightarrow} Q_1 \stackrel{\partial_1'}
{\longrightarrow} Q_0 \stackrel{\epsilon'}{\dashrightarrow} M
\dashrightarrow 0 $\end{center}

\noindent with the $P_i$ and $Q_i$ all projective modules.  Our
main result is:

\begin{thm} \hspace{2mm} The complexes:

\[ P_n \oplus S_n \stackrel{\partial_n \oplus 0}{\longrightarrow}
P_{n-1} \stackrel{\partial_{n-1}}{\longrightarrow} \cdots
\stackrel{\partial_2}{\longrightarrow} P_1
\stackrel{\partial_1}{\longrightarrow} P_0 \eqno(1)
\]
\noindent and \[ Q_n \oplus T_n \stackrel{\partial_n'\oplus
0}{\longrightarrow} Q_{n-1}
\stackrel{\partial_{n-1}'}{\longrightarrow} \cdots
\stackrel{\partial_2'}{\longrightarrow}Q_1
\stackrel{\partial_1'}{\longrightarrow} Q_0 \eqno(2) \]

\bigskip
\noindent are chain homotopy equivalent, where the projective
modules $T_i$, $S_i$ are defined inductively by:

\bigskip
\noindent$T_0 \cong P_0$,  $S_0 \cong Q_0$ and for $i=1,\dots, n$:
$T_i  \cong  S_{i-1} \oplus P_i$

\hspace{56mm}$S_i  \cong  T_{i-1} \oplus Q_i$

\end{thm}

Given maps $f:A \to C$, $g:B \to C$ the notation $f \oplus g$ will
always be used to denote the map $f\oplus g:A \oplus B \to C$
given by $f\oplus g: (a,b) \mapsto f(a)+g(b)$.

\section{Construction of chain homotopy equivalence}

For each $i\in 1,\dots, n$ we have natural inclusions of summands:
\[\iota_i:P_i \to T_i (\cong  S_{i-1} \oplus P_i) \qquad\qquad
\iota_i':Q_i \to S_i  (\cong  T_{i-1} \oplus Q_i)\] Let
$\iota_0:P_0 \to T_0$ and $\iota_0':Q_0 \to S_0$ both be the
identity map.

\bigskip
 For $i=1,\dots, n$, we define $\delta_i:T_i (\cong
P_i \oplus S_{i-1}) \to T_{i-1} \oplus S_{i-1}$

\hspace{33mm} and $\delta_i':S_i (\cong Q_i \oplus T_{i-1}) \to
S_{i-1}\oplus T_{i-1}$

\noindent by \[\delta_i=\left(
\begin{array}{cc}
\iota_{i-1}\partial_i&0\\ 0&1 \end{array} \right)\qquad\qquad
\delta_i'=\left(
\begin{array}{cc}
\iota_{i-1}'\partial_i'&0\\ 0&1 \end{array} \right)\]

\bigskip
For $r=0,\dots, n-1$, let $\mathcal{C}_r$ denote the chain
complex:

 \[
    P_n \oplus S_n \stackrel{\partial_n \oplus 0}{\longrightarrow}
     \cdots
     \stackrel{\partial_{r+2}}{\longrightarrow}
     P_{r+1}
    \stackrel{\iota_{r}\partial_{r+1}}{\longrightarrow}T_{r}
    \stackrel{\delta_{r}}{\longrightarrow} T_{r-1} \oplus S_{r-1}
    \stackrel{\delta_{r-1}\oplus 0}{\longrightarrow}  \cdots
    \stackrel{\delta_1 \oplus 0}{\longrightarrow} T_{0} \oplus S_0
    \]

Also let $\mathcal{C}_n$ denote the chain complex:

\[
T_n \oplus S_n \stackrel{\delta_n \oplus 0}{\longrightarrow}
T_{n-1} \oplus S_{n-1} \stackrel{\delta_{n-1} \oplus 0}
{\longrightarrow} \cdots \cdots  \stackrel{\delta_2 \oplus
0}{\longrightarrow} T_1 \oplus S_1 \stackrel{\delta_1 \oplus
0}{\longrightarrow} T_{0} \oplus S_0
\]

  Clearly $\mathcal{C}_0$ is the chain complex (1).
For $r=0, \dots n-1$, the chain complex $\mathcal{C}_{r+1}$ is
obtained from $\mathcal{C}_{r}$ by replacing:
\[\stackrel{\partial_{r+2}}{\longrightarrow} P_{r+1}
\stackrel{\iota_{r}\partial_{r+1}}{\longrightarrow}T_{r}
\stackrel{\delta_r}{\longrightarrow}\] with
\[\stackrel{\iota_{r+1}\partial_{r+2}}{\longrightarrow} P_{r+1}
\oplus S_r \stackrel{\delta_{r+1}}{\longrightarrow} T_{r} \oplus
S_{r}\stackrel{\delta_r \oplus 0}{\longrightarrow}\quad\,\]

This is a simple homotopy equivalence so $\mathcal{C}_{r+1}$ is
chain homotopy equivalent to $\mathcal{C}_{r}$.

\bigskip
Similarly, for $r=0,\dots, n-1$, let $\mathcal{D}_r$ denote the
chain complex:

 \[ Q_n \oplus T_n \stackrel{\partial_n' \oplus 0}{\longrightarrow}
    \cdots \stackrel{\partial_{r+2}'}{\longrightarrow} Q_{r+1}
    \stackrel{\iota_{r}'\partial_{r+1}'}{\longrightarrow}S_{r}
    \stackrel{\delta_{r}'}{\longrightarrow} S_{r-1} \oplus T_{r-1}
    \stackrel{\delta_{r-1}'\oplus 0}{\longrightarrow}  \cdots \stackrel{\delta_1' \oplus 0}{\longrightarrow} S_{0}
    \oplus T_0
    \]

Again let $\mathcal{D}_n$ denote the chain complex:

\[
S_n \oplus T_n \stackrel{\delta_n' \oplus 0}{\longrightarrow}
S_{n-1} \oplus T_{n-1} \stackrel{\delta_{n-1}' \oplus 0}
{\longrightarrow} \cdots  \stackrel{\delta_2' \oplus
0}{\longrightarrow} S_1 \oplus T_1 \stackrel{\delta_1' \oplus
0}{\longrightarrow} S_{0} \oplus T_0
\]

  Clearly $\mathcal{D}_0$ is the chain complex (2). As
before, for $r=0, \dots n-1$, the chain complex
$\mathcal{D}_{r+1}$ is chain homotopy equivalent to
$\mathcal{D}_{r}$.

\bigskip
  We have (1) chain homotopy equivalent to
$\mathcal{C}_{n}$ and (2) chain homotopy equivalent to
$\mathcal{D}_{n}$.  We complete the proof of the theorem by
showing that $\mathcal{C}_{n}$ is chain isomorphic to
$\mathcal{D}_{n}$.

\begin{lemma}  There exist inverse pairs of maps $h_i$, $k_i$ making
the following diagram commute:

\begin{eqnarray*} T_n \oplus S_n
\stackrel{\delta_n \oplus 0}{\longrightarrow} T_{n-1} \oplus
S_{n-1} \stackrel{\delta_{n-1} \oplus 0} {\longrightarrow} \cdots
\cdots  \stackrel{\delta_2 \oplus 0}{\longrightarrow} T_1 \oplus
S_1 \stackrel{\delta_1 \oplus 0}{\longrightarrow} T_{0} \oplus S_0
\stackrel{\epsilon \oplus 0}
{\dashrightarrow} M \dashrightarrow 0\\
\downarrow h_n \,\,\,\,\,\qquad\qquad \downarrow h_{n-1}\quad
\,\,\qquad\qquad \qquad\qquad \downarrow h_1 \,\,\,\quad\quad\quad
\downarrow h_0 \,\,\quad\quad\downarrow 1\qquad
\\
S_n \oplus T_n \stackrel{\delta_n' \oplus 0}{\longrightarrow}
S_{n-1} \oplus T_{n-1} \stackrel{\delta_{n-1}' \oplus 0}
{\longrightarrow} \cdots  \stackrel{\delta_2' \oplus
0}{\longrightarrow} S_1 \oplus T_1 \stackrel{\delta_1' \oplus
0}{\longrightarrow} S_{0} \oplus T_0 \stackrel{\epsilon' \oplus 0}
{\dashrightarrow} M \dashrightarrow 0\\
\downarrow k_n \,\,\,\,\,\qquad\qquad \downarrow k_{n-1}\quad
\,\,\qquad\qquad \qquad\qquad \downarrow k_1 \,\,\,\quad\quad\quad
\downarrow k_0 \,\,\quad\quad\downarrow 1\qquad
\\
T_n \oplus S_n \stackrel{\delta_n \oplus 0}{\longrightarrow} T_{n-1}
\oplus S_{n-1} \stackrel{\delta_{n-1} \oplus 0} {\longrightarrow}
\cdots  \stackrel{\delta_2 \oplus 0}{\longrightarrow} T_1 \oplus S_1
\stackrel{\delta_1 \oplus 0}{\longrightarrow} T_{0} \oplus S_0
\stackrel{\epsilon \oplus 0} {\dashrightarrow} M \dashrightarrow 0
\end{eqnarray*}

\end{lemma}

\bigskip
 Proof:  As $T_0$, $S_0$ are projective, we may
pick $f_0$, $g_0$ so that the following diagrams commute:
\begin{eqnarray*}
T_0 \stackrel{\epsilon}{\longrightarrow} M\,\,   \qquad \qquad
T_0 \stackrel{\epsilon}{\longrightarrow} M\,\, \nonumber\\
\downarrow f_0 \,\, \quad \downarrow 1             \qquad \qquad
\uparrow g_0 \,\, \quad \uparrow 1
\nonumber\\
S_0 \stackrel{\epsilon'}{\longrightarrow} M\,\,  \qquad \qquad S_0
\stackrel{\epsilon'}{\longrightarrow} M\,\,
\end{eqnarray*}
\hfill (3)

Define $h_0:T_0 \oplus S_0 \to S_0 \oplus T_0$ and $k_0:S_0 \oplus
T_0 \to T_0 \oplus S_0$ by:

\[
h_0= \left(\begin{array}{cc} f_0&1-f_0g_0\\ 1& -g_0
\end{array}\right)
\qquad \qquad k_0= \left(\begin{array}{cc} g_0&1-g_0f_0\\ 1& -f_0
\end{array}\right)
\]

\bigskip
\noindent Direct calculation shows that $h_0k_0 = 1$ and $k_0h_0 =
1$.

\bigskip
 Also from commutativity of (3), we deduce:

\[
(\epsilon' \quad 0)\left(\begin{array}{cc} f_0&1-f_0g_0\\ 1& -g_0
\end{array}\right)
= (\epsilon'f_0 \quad \epsilon'(1-f_0g_0))=(\epsilon \quad 0)
\]
and
\[
(\epsilon \quad 0)\left(\begin{array}{cc} g_0&1-g_0f_0\\ 1& -f_0
\end{array}\right)
= (\epsilon g_0 \quad \epsilon(1-g_0f_0))=(\epsilon' \quad 0)
\]

\bigskip

Hence the following diagrams commute:

\begin{eqnarray*}
T_0 \oplus S_0 \stackrel{\epsilon \oplus 0}{\longrightarrow} M\,\,
\qquad \qquad       T_0 \oplus S_0
\stackrel{\epsilon \oplus 0}{\longrightarrow} M\,\,\\
\downarrow h_0 \,\,\quad\,\, \quad \downarrow 1 \quad \qquad
\qquad \uparrow k_0 \,\, \quad\,\, \quad
\uparrow 1\\
S_0 \oplus T_0 \stackrel{\epsilon'\oplus 0}{\longrightarrow} M\,\,
\qquad \qquad       S_0 \oplus T_0 \stackrel{\epsilon' \oplus
0}{\longrightarrow} M\,\,\\
\end{eqnarray*}

 Now suppose that for some $0<i \leq n$, we have defined
$h_j:T_j \oplus S_j \to S_j \oplus T_j$ and $k_j:S_j \oplus T_j
\to T_j \oplus S_j$ for $j=0,\dots,i-1$, so that for each $j$, we
have $h_jk_j=1$ and $k_jh_j=1$.  We proceed by induction.

As before, pick $f_i$, $g_i$ so that the following diagrams
commute:

\begin{eqnarray*}
T_i \stackrel{\delta_i}{\longrightarrow} T_{i-1} \oplus
S_{i-1}\,\, \qquad \qquad    T_i
\stackrel{\delta_i}{\longrightarrow}
T_{i-1} \oplus S_{i-1}\,\, \nonumber\\
\downarrow f_i \,\quad \quad \, \quad \downarrow h_{i-1}
\,\,\,\qquad \qquad \uparrow g_i \,\, \quad \quad \quad
\uparrow k_{i-1} \,\,\,\,\nonumber\\
S_i \stackrel{\delta_{i}'}{\longrightarrow} S_{i-1} \oplus T_{i-1}
\,\, \qquad \qquad   S_i \stackrel{\delta_{i}'}{\longrightarrow}
S_{i-1} \oplus T_{i-1} \,\,
\end{eqnarray*}
\hfill(4)

Define $h_i:T_i \oplus S_i \to S_i \oplus T_i$ and $k_i:S_i \oplus
T_i \to T_i \oplus S_i$ by:

\[
h_i= \left(\begin{array}{cc} f_i&1-f_ig_i\\ 1& -g_i
\end{array}\right) \qquad \qquad
k_i= \left(\begin{array}{cc} g_i&1-g_if_i\\ 1& -f_i
\end{array}\right)
\]

\bigskip
\noindent  Direct calculation shows that $h_ik_i = 1$ and $k_ih_i
= 1$.
\bigskip

 Recall $h_{i-1}k_{i-1}=1$ and $k_{i-1}h_{i-1}=1$.
From commutativity of (4) we deduce:

\[
(\delta_i' \quad 0)\left(\begin{array}{cc} f_i&1-f_ig_i\\ 1& -g_i
\end{array}\right)
= (\delta_i'f_i \quad \delta_i'(1-f_ig_i))=h_{i-1}(\delta_i \quad
0)
\]

\noindent and

\[
(\delta_i \quad 0)\left(\begin{array}{cc} g_i&1-g_if_i\\ 1& -f_i
\end{array}\right)
= (\delta_i g_i \quad \delta_i(1-g_if_i))=k_{i-1}(\delta_i' \quad
0)
\]

Hence the following diagrams commute:

\begin{eqnarray*}
T_i \oplus S_i \stackrel{\delta_i \oplus 0}{\longrightarrow}
T_{i-1} \oplus S_{i-1}\,\, \qquad \qquad    T_i \oplus S_i
\stackrel{\delta_i \oplus 0}{\longrightarrow}
T_{i-1} \oplus S_{i-1}\,\,\\
\downarrow h_i \,\quad \quad \, \quad \quad \downarrow h_{i-1}
\,\,\,\qquad \quad \qquad \uparrow k_i \,\, \quad \quad \quad
\quad
\uparrow k_{i-1}\,\,\,\,\\
S_i \oplus T_i \stackrel{\delta_{i}' \oplus 0}{\longrightarrow}
S_{i-1} \oplus T_{i-1} \,\, \qquad \qquad   S_i \oplus
T_i\stackrel{\delta_{i}' \oplus 0}{\longrightarrow} S_{i-1} \oplus
T_{i-1} \,\,
\end{eqnarray*}

\bigskip
So we may construct the $h_i$, $k_i$ as required. \hfill $\Box$
\,\,

\bigskip
We know the $h_i$, $i=0, \dots,n$ constitute a chain map
$h:\mathcal{C}_n \to \mathcal{D}_n$.  Also the $k_i$ constitute a
chain map $k:\mathcal{D}_n \to \mathcal{C}_n$.  As $h$ and $k$ are
mutually inverse we have that $\mathcal{C}_n$ and $\mathcal{D}_n$
are chain isomorphic.  Hence (1) and (2) are chain homotopy
equivalent as required.

\section{Injective Resolutions}
Finally we note that dual arguments may be used in the same way to
prove the dual result:

\begin{thm}  Let $(I_r, \partial_r)$ and
$(J_r, \partial_r')$ be injective resolutions for some module $M$,
truncated after the $n^{\rm th}$ terms (so $M \cong {\rm
Ker}(\partial_0: I_0 \to I_1) \cong {\rm Ker}(\partial_0': J_0 \to
J_1$)). Then stabilizing the final modules, $I_n$ and $J_n$, with
the appropriate injective modules results in chain homotopy
equivalent complexes. \end{thm}

email: {\verb|wajid@mannan.info| }

Address: Mathematics Department, UCL, Gower Street - London - WC1E 6BT

\end{document}